\documentstyle[amsfonts,12pt]{article}
\newtheorem{theorem}{\sc Theorem}
\newtheorem{proposition}[theorem]{\sc Proposition}
\newtheorem{lemma}[theorem]{\sc Lemma}
\newtheorem{corollary}[theorem]{\sc Corollary}

\newcommand{\hao}{\ha_1}
\newcommand{\goo}{\go_1}
\newcommand{\restriction}{{|}} 

\newcommand{\proof}{\noindent {\sc Proof. }} 
 
\def\rest{\mathord{\restriction}}

\newcommand{\open}{\Bbb}
 
\newlength{\labparwidth}
\renewcommand{\hom}{\mbox{\rm Hom}}

\newcommand{\rge}{\mbox{\rm rge}}
\newcommand{\ext}{\mbox{\rm Ext}}

\newcommand{\se}{\subseteq}
\newcommand{\set}[2]{\{#1 \colon #2\}} 
\newcommand{\fin}{$\Box$\par\medskip} 

\renewcommand{\to}{\rightarrow}
 
\newcommand{\ga}{\alpha}
\newcommand{\gb}{\beta}
\newcommand{\grg}{\gamma}
\newcommand{\gd}{\delta}
\newcommand{\gre}{\varepsilon}
\newcommand{\gz}{\zeta}
\newcommand{\gh}{\eta}
\newcommand{\gth}{\theta}

\newcommand{\gl}{\lambda}

\newcommand{\gn}{\nu}

\newcommand{\gs}{\sigma}

\newcommand{\gf}{\varphi}

\newcommand{\go}{\omega}
 
\newcommand{\gG}{\Gamma}

\newcommand{\gF}{\Phi}

\newcommand{\ha}{\aleph}
 
\newcommand{\oZ}{{\open Z}}

\title{Uniformization and the Diversity of Whitehead Groups}
\author{P. C. Eklof \\University of California, Irvine \and A. H.
Mekler \\Simon Fraser University\thanks{Research partially supported
by NSERC grant \#9848} \and S. Shelah \\Hebrew University
and \\Rutgers University\thanks{Research partially supported by the
BSF. The authors thank Rutgers
University for its support. Publication \#441}}

\begin{document}
\maketitle

\section*{Introduction}

The connections between Whitehead groups and uniformization
properties were investigated by the third author in  \cite{Sh98}. In
particular it was essentially shown there that there is a non-free
Whitehead  (respectively, $\hao$-coseparable) group of cardinality
 $\hao$ if and only if there is a
ladder system on a stationary subset of $\goo$ which satisfies
$2$-uniformization (respectively, $\go$-uniformization).
 (See also \cite[\S XII.3]{EM}; definitions are
reviewed below.) These techniques allowed also the proof of various
independence and consistency results about Whitehead groups, for
example that it is consistent that there is a non-free Whitehead
group of cardinality $\hao$ but no non-free $\hao$-coseparable group
(cf.\ \cite[XII.3.18]{EM}).

However, some natural questions remained open, among them the
following two, which are stated as problems at the end of
\cite[p. 454]{EM}. 
\begin{itemize}
\item  Is it consistent that the class of W-groups of cardinality $\aleph _1$ is 
exactly the class of strongly $\aleph _1$-free groups of cardinality 
$\aleph _1$?
\item  If every strongly $\aleph _1$-free group of cardinality $\aleph _1$ is a 
W-group, are they also all $\aleph _1$-coseparable? 
\end{itemize}
In this paper we use the techniques of uniformization to answer the
first question in the negative and give a partial affirmative answer
to the second question. (The third author claims a full affirmative
solution to the second question, but it is too complicated to give here.)

More precisely, we have the following two theorems of ZFC.

\begin{theorem}
\label{sumup}
The following are equivalent:

\noindent 
(a) There is an $\aleph _1$-separable Whitehead group A of
cardinality $\aleph _1$ with $\Gamma (A) = 1$.

\noindent 
(b) There is a strongly $\aleph _1$-free Whitehead group A
of cardinality $\aleph _1$ with $\Gamma (A) = 1$.

\noindent
(c) There is a Whitehead group A of cardinality $\aleph _1$ with
$\Gamma (A) = 1$.

\noindent
(d) There is a  Whitehead group of cardinality 
$\aleph _1$ which is {\em not} strongly $\aleph _1$-free.

\noindent
(e) There is a ladder system on $\lim (\omega _1)$ 
which satisfies 2-uniformization.

\end{theorem}
 
The new part of this result is the proof of (d) from (c); this gives
a negative answer to the first question. Given the history of
independence results regarding Whitehead groups, it is remarkable
that the answer to this question is
negative.
\footnote{Alan: this sentence was added to the version submitted. I
think you didn't get the last version, of Feb. 28, with this and a few
other minor changes. I managed to find it. (Problem was I had not kept
up the revision history at the top, so that confused me.)}
The partial
answer to the second question is contained in the following.

\begin{theorem}
\label{sumup3}
Consider the following hypotheses.

(1) Every strongly $\ha_1$-free group of cardinality $\ha_1$ is
 $\ha_1$-coseparable.

(2)  Every strongly $\ha_1$-free group of cardinality $\ha_1$ is
a Whitehead group.

(3) Every ladder system on a stationary subset of $\go_1$ satisfies
$2$-uniformization.

(4) Every ladder system on a stationary subset of $\go_1$ satisfies
$\go$-uniformization.

(5) There is  a strongly $\ha_1$-free group of cardinality $\ha_1$
which  is
 $\ha_1$-coseparable but not free.

\noindent
Then (1) $\Rightarrow$ (2) $\Rightarrow$ (3) $\Leftrightarrow$ (4)
$\Rightarrow$ (5).

\end{theorem}

The new parts of this theorem are the proofs of (3) from (2) and (4)
from (3). We consider
 the implication from (2) to
(4) strong evidence for an affirmative answer to the
second question; what is lacking for a complete answer
 is a proof of (1) from (4). The implication from (2) to (5) is
also new and of interest.

The last two sections of this paper contain some other results about
uniformization, which may be of independent interest.

\section*{Preliminaries}

Let us review some basic notation and terminology. See \cite{EM} for
further information; throughout the paper we will usually cite
\cite{EM} for results we need, rather than the original source.

We will always be dealing with abelian groups or $\oZ$-modules;  we shall 
simply  say ``group''. A group $A$ is said to be a {\em Whitehead group} if
$\ext(A, \oZ) = 0$; it is said to be {\em $\hao$-coseparable} if
 $\ext(A, \oZ^{(\go)}) = 0$. 

A  group $A$ of arbitrary cardinality is 
{\em $\hao$-free} if and only if every countable subgroup of $A$ is 
free; $A$ is {\em strongly $\hao$-free} if and only if every countable subset
 is contained in a free subgroup $B$ such that $A/B$ is $\hao $-free.
$A$ is {\em $\hao$-separable}
 if and only if every countable subset
 is contained in a free subgroup $B$ such that $B$ is a direct summand
of $A$.
\footnote{referee: ``Since you bothered to give the definition of
the basic concepts like ``$\hao$-free'' etc., include the definition of
``$\hao$-separable''.''}

Chase \cite{C} showed that CH implies that every Whitehead group 
is  strongly $\hao$-free.
In the third author's original paper, \cite{Sh44}, on the
independence of the Whitehead Problem, a larger class of groups than
the strongly $\hao$-free groups plays a key role, namely the groups
which the first author (\cite{E}) later named the {\em Shelah groups}.
These are the $\hao$-free groups $A$ such that for every countable
subgroup $B$ there is a countable subgroup $B'\supseteq B$  such that
for any countable $C$ satisfying  $C \cap  B' = B$,
 $C/B$ is free. In \cite{Sh44} it is proved consistent --- in fact a
consequence of Martin's Axiom plus $\neg$CH --- that every Shelah
group of cardinality $\hao$ is $\hao$-coseparable. Later, in
\cite{Sh105} it was proved consistent
 --- in fact, again a consequence of Martin's Axiom plus
$\neg$CH
 ---  that the Whitehead groups of cardinality $\hao$ are the same as
the $\hao$-coseparable groups and are precisely the Shelah groups. The
first author emphasized the strongly $\hao$-free groups in his
expository accounts of this work (e.g. in \cite{E',E}),
 as a class of groups more  familiar to 
algebraists, and raised the first question cited above. The answer to
that question now given here now shows, definitively,
 that the larger class of Shelah groups is the `right one' to consider
for the Whitehead Problem.

Notions of uniformization (in our sense) were first defined in
\cite{Sh65}.   Let $S$ be a subset
of $\lim (\omega _1)$. If  $\delta  \in  S$, a {\it ladder on $\delta
$} is a function $\eta _\delta \colon \omega  \rightarrow \delta $
which is strictly increasing and has range cofinal in $\delta $. A
{\it ladder system on $S$} is an indexed family  $\eta  = \{\eta
_\delta \colon \delta  \in  S\}$  such that each $\eta _\delta $ is a
ladder on $\delta $.  For a cardinal  $\lambda  \geq  2$, a $\lambda
${\it -coloring} of a ladder system $\eta $ on $S$ is a family  $c =
\{c_\delta \colon \delta  \in  S\}$  such that  $c_\delta \colon
\omega  \rightarrow  \lambda $.  A {\it uniformization} of a coloring
$c$ of a ladder system $\eta $ on $S$ is a pair $\langle f,
f^*\rangle $ where
$f\colon \omega _1 \rightarrow  \lambda $, $f^*\colon S \rightarrow
\omega $  and for all  $\delta \in  S$  and all  $n \geq  f^*(\delta
)$, $f(\eta _\delta (n)) = c_\delta (n)$. If such a pair exists, we
say that $c$ can be uniformized. In order for the pair to exist it is
enough to have either member of the pair; i.e., either $f$ so that for
all  $\delta  \in  S$, $f(\eta _\delta (n)) = c_\delta (n)$, for all
but finitely many $n$, or $f^*$ so that for all  $\delta $, $\alpha
\in  S$, if  $n \geq  f^*(\delta )$, $m \geq f^*(\alpha )$  and  $\eta
_\delta (n) = \eta _\alpha (m)$, then  $c_\delta (n) = c_\alpha (m)$.
We say that $(\eta, \lambda)${\it -uniformization holds} or that {\it
$\eta$ satisfies $\gl$-uniformization } if every $\lambda $-coloring
of $\eta $ can be uniformized. We will generalize these (by now,
standard) notions in the next section.

If $A$ is an $\hao$-free group of cardinality $\hao$, then (we
define) $\gG(A) =
1$ if and only if $A$ is the union of a continuous chain of countable
subgroups
$$ A = \bigcup_{\ga < \goo} A_\ga$$
such that for all $\ga \in \lim(\goo)$, $A_{\ga + 1}/A_\ga$ is not free.
If $A$ is not strongly $\hao$-free, then $\gG(A) = 1$, but the
converse is false.

\begin{lemma}

\label{XII.3}

If there is a Whitehead
group $A$ of cardinality $\hao$ with $\gG(A) = 1$, then there is a
ladder system on $\lim(\goo)$ which satisfies 2-uniformization.

\end{lemma}

\proof
We assume familiarity with \cite[\S XII.3]{EM} and sketch the
modifications to the proof of Theorem XII.3.1 that are needed.
 In the proofs of Lemma XII.3.16 and Theorem XII.3.1,
$\lim(\goo)$ is partitioned into countably many sets $E_n$; to each of
these is associated $\gF^n = \{\gf_\ga: \ga \in E_n\}$, which is a
family with $2$-uniformization. As defined there, the range of
the $\gf_\ga$ is not a set of ordinals, but it is easy to see that,
by a coding argument, we can assume that the range of $\gf_\ga$ is
contained in $\ga$ and, furthermore, that if $\ga \in E_i$ and $\gb \in E_j$,
then the ranges of $\gf_\ga$ and $\gf_\gb$ are disjoint. Finally, if necessary,
one modifies each $\gf_\ga$ so that it is a ladder on $\ga$ (say by
using a bijection from $\goo \times \goo$ to $\goo$). This produces a
ladder system on $\lim(\goo)$ which has $2$-uniformization since the
uniformizations of the original $\gF^n$ fit together to give a
uniformization of the ladder system.
\fin

This proof obviously generalizes to prove that if there is a Whitehead
group $A$ of cardinality $\hao$ with $\gG(A) = \tilde{S}$, then there is a
ladder system on $S$ which satisfies 2-uniformization.

If $\ga < \gb$ are ordinals, denote by $(\ga, \gb)$ the open interval
of ordinals between $\ga$ and $\gb$, i.e., the set $\{\grg \colon \ga
< \grg < \gb \}$. Similarly we define the half open interval $[\ga,
\gb)$, etc. We will use $\langle \ga, \gb \rangle$ to denote the ordered
pair of ordinals.



\section{The First Question}

It is consistent that every strongly $\hao$-free group of cardinality
$\hao$ is Whitehead (cf.\ \cite [XII.1.12]{EM}) and it is consistent
that there are non-free Whitehead groups of cardinality $\hao$ and
every Whitehead group of cardinality $\hao$ is strongly $\hao$-free
(cf.\ \cite[XII.1.9]{EM}), but here we show that it's not consistent
that the Whitehead groups of cardinality $\hao$ are precisely the
strongly $\hao$-free groups.

  If $\alpha  \in  [\omega, \omega _1)$ and $\alpha  = \delta
+ n$, where $\delta $ is a limit ordinal and $n \in  \omega $, a {\it
ladder on $\alpha $} is defined to be a ladder on $\gd$.
 Thus, for example, a ladder on $\omega  + 1$ is a strictly
increasing $\omega$-sequence approaching $\omega $. If $S \subseteq
[\omega, \omega _1)$, a {\it ladder system on $S$} is an indexed
family $\eta  = \langle \eta _\alpha \colon \alpha  \in S\rangle $ such
that each $\eta _\alpha$ is a ladder on $\alpha$. 

Whenever we write an ordinal as $\delta  + n$ we mean that $\delta
\in \lim(\omega _1)$ and $n \in  \omega$. We will always assume in
what follows that if $\delta  + n \in  S$, then $\delta  \in  S.$ 

Suppose that $H$ is an indexed family $\langle h_\alpha \colon \alpha \in
S\rangle$ where each $h_\alpha$ is a function: $\omega \rightarrow
\omega$. If $\eta = \langle \eta _\alpha \colon \alpha \in S\rangle $ is
a ladder system on $S$, an $H${\it -coloring of $\eta $} is an
indexed family $c = \langle c_\alpha \colon \alpha \in S\rangle $ such
that for all $\alpha $, $c_\alpha \colon \omega \rightarrow \omega $ and
that for all $n \in \omega $, $c_\alpha (n) < h_\alpha (n)$. We say
that $(\eta, H)${\em -uniformization holds} (or $\eta $ {\it
satisfies $H$-uniformization}) if whenever $c$ is an $H$-coloring,
there is a pair $\langle f, f^\ast\rangle$ such that $f\colon \omega_1 \rightarrow
\omega$, $f^\ast \colon S \rightarrow \omega $, and for all $\alpha \in
S$, $f(\eta_\alpha (n)) = c_\alpha (n)$ whenever $n \geq f^\ast
(\alpha )$. We say that $(\eta, \gl)$-uniformization holds if each
$h_\alpha \in H$ is the constant function $\gl$; this agrees with
the previous definition. 

\medskip 
A ladder system
$\eta = \langle \eta _\alpha \colon \alpha \in S\rangle $ is said to be
{\em tree-like\/} if for all $\alpha $, $\beta \in S$, if $\eta _\alpha (n) =
\eta _\beta (m)$, then $n = m$ and $\eta _\alpha (k) = \eta _\beta
(k)$ for all $k \leq n$. Let $F$ be a function from $S$ to $\omega$;
say that $\eta$ is {\em strongly tree-like w.r.t.\/} $F$ if $\eta$ is
tree-like and in addition, whenever $\eta_\alpha(n) = \eta_\beta(m)$
for some $\alpha, \beta \in S$ and $n, m \in \omega$, then
$F(\alpha) = F(\beta)$.

\begin{lemma} 
\label{tree}
Suppose that there is a ladder system $\zeta  = \langle \zeta _\ga
\colon \ga  \in  S\rangle $ on $S \supseteq  \lim (\omega _1)$ such
that $(\zeta, H)$-uniformization holds. Given a function $F\colon S
\rightarrow  \omega$, there is a ladder system $\eta  = \langle \eta
_\alpha \colon \alpha  \in  S\rangle$ such that $\eta$ is strongly
tree-like w.r.t.\  $F$  and $(\eta, H)$-uniformization holds. 
\end{lemma}

\proof
 Choose  a one-one onto function $\gth$ from $\go
\times {}^{<\go}\go_1$ to $\go_1$ 
with the property that for all limit $\gd$, $\gth[\go
\times {}^{<\go}\gd] = \gd$ and for all $k \in \go$, if $t$ is a sequence which extends
$s$ then $\gth(k, s) < \gth(k, t)$. For each $\ga$, let $\eta_\ga(n) =
 \gth(\langle F(\ga), \langle \zeta_\ga(m)
\colon m \leq n\rangle \rangle)$.
 Since $\gth(k, s) < \gth(k, t)$, $\eta_\ga$ is strictly
increasing.
If we can show that each $\gz_\ga$ is a ladder on $\ga$, then we will
be done since, by construction, it is strongly tree-like w.r.t. $F$.
 Observe that because $\gth\colon \go\times{}^{<\go}\gd \to \gd$
is one-one and onto for limit $\gd$, if $\mu$ is a limit
ordinal $ \leq \zeta_\ga(n)$, then $\gth(\langle F(\ga), \langle \zeta_\ga(m)
\colon m \leq n\rangle \rangle) \geq \mu$. Consider now $\ga = \gd + n$.
Note that
$\eta_\ga$ has range contained in $\gd$. If
$\gd$ is a limit of limit ordinals then, by the observation, the range
of $\eta_\ga$ is cofinal in $\gd$ since the range of $\zeta_\ga$ is
cofinal. If $\gd = \grg + \go$ then there is some $k$ so that
$\zeta_\ga(k) \geq \grg$. Then for all $ m\geq k$, $\grg \leq\eta_\ga(m)
< \gd = \grg + \go$. So $\eta_\ga$ is cofinal in $\gd$. 
\fin

\begin{lemma} 
\label{prep}

Suppose that there is a ladder system $\zeta = \langle \zeta _\ga
\colon \ga \in S\rangle $ where $S \supseteq \lim(\omega _1)\rangle $
such that $(\zeta, 2)$-uniformization holds. Given $H = \langle
h_\alpha \colon \alpha \in [\omega, \omega _1)\rangle $ where each
$h_\alpha \colon \omega \rightarrow \omega$ and given a function $F\colon
[\omega, \omega _1) \rightarrow \omega $, there is a ladder system
$\eta = \langle \eta _\alpha \colon \alpha \in [\omega $,
$\omega_1)\rangle$ such that $(\eta, H)$-uniformization holds and
$\eta$ is strongly tree-like w.r.t.\ $F$. 
\end{lemma} 

\proof
We shall give the proof as a series of reductions. First of all, by
\cite[XII.3.2]{EM}, $(\zeta, 3)$-uniformization holds. Next, we claim
that we can assume that $\zeta$ is a ladder system on $[\omega,
\omega _1)$. Write $\omega$ as the disjoint union of $\aleph_0$
disjoint infinite sets $Y_n\ (n \in \omega)$, and for each $n$ let
$\theta_n\colon \omega \rightarrow Y_n$ enumerate $Y_n$ in increasing
order. For each $\delta \in \lim (\omega _1)$ and $n \in \omega $,
define $\zeta'_{\delta +n} = \zeta _\delta \circ \theta _n$. Then it
is easy to see that $\langle \zeta '_\alpha \colon \alpha \in [\omega,
\omega _1)\rangle$ satisfies 3-uniformization. 

So we will now assume that $S = [\omega, \omega _1)$. By
Lemma~\ref{tree}, we can assume that $\zeta$ is tree-like.  For each
$\alpha \in S$ define $\psi_\alpha \colon \omega \rightarrow \omega $
by $\psi _\alpha (n) = \Sigma _{j\leq n} h_\alpha (j)$; so  $\psi
_\alpha (n) - \psi _\alpha (n-1) = h_\alpha (n)$ for all $n \in \omega
$ (where $\psi _\alpha (-1) = 0)$. Define $\zeta'_\alpha = \zeta
_\alpha \circ \psi _\alpha $. Now we claim that we can assume that
$\zeta'$ satisfies $H$-uniformization.
Suppose that $c' = \langle c'_\alpha \colon \alpha \in S\rangle $ is
an $H$-coloring of $\zeta'$. Define a 3-coloring $c$ of $\zeta $ as
follows. Let   $c_\alpha (0) = 2$, and for  each $n \in \omega $, and $\alpha \in S$, let $c_\alpha
(\psi _\alpha (n)) = 2$.  Define $c_\alpha (\psi
_\alpha (n-1) + k) = 0$ for $1 \leq k \leq c_\alpha '(n)$, and
$c_\alpha (\psi _\alpha (n-1) + k) = 1$ for $c_\alpha '(n) < k <
h_\alpha (n).$ 

As an example, suppose $h_\alpha (0) = 5$, $h_\alpha (1) = 4$,
$h_\alpha (2) = 5$ and $h_\alpha (3) = 6$. Then $\psi _\alpha (0) =
5$, $\psi _\alpha (1) = 9$, $\psi _\alpha (2) = 14$ and $\psi _\alpha
(3) = 20$. If $c_\alpha'(0) = 4$, $c_\alpha '(1) = 1$, $c_\alpha
'(2) = 3$ and $c_\alpha'(3) = 0$, then the values of $c_\alpha (n)$
for $0 \leq n \leq 20$ are: 
$$ 2,0,0,0,0,2,0,1,1,2,0,0,0,1,2,1,1,1,1,1,2. $$
(The blocks of 0's between 2's code the values of $c_\alpha'$.)

 Given
$\langle f, f^*\rangle$ which uniformizes $c$, define $\langle f', {f'}^*\rangle$ as follows.
Let $n \geq {f'}^*(\alpha)$, if and only if $\psi _\alpha (n-1)  \geq
f^*(\alpha)$. We need to choose $f'$ so that
$f'(\nu) = c'_\ga(m)$ if $\nu = \gz'_\ga(m)$ and $m \geq f^{'*}(\ga)$.
To see that there is such an $f'$, suppose $\gb$ and $k$ are such
that also
$\nu = \gz'_\gb(k)$ where $k \geq f^{'*}(\gb)$. 
Since $\langle \gz_i\colon i \in [\go,
\go_1)\rangle$ is tree-like we have that $\gz_\ga\rest(\psi_\ga(m) +
1) = \gz_\gb\rest(\psi_\ga(m) + 1)$ and  $\psi_\ga(m) =
\psi_\gb(k)$. For definiteness assume that $\psi_\ga(m-1) \leq
\psi_\gb(k-1)$. Then since $\psi_\ga(m-1) \geq f^*(\ga)$ and
$\psi_\gb(k-1) \geq f^*(\gb)$, we have that $c_\ga(r) = c_\gb(r)$ for
all $r$ such that $\psi_\gb(k-1) \leq r \leq \psi_\gb(k)$. By
the coding we know that $\psi_\ga(m-1)$ is the greatest natural
number, $s$, less than $\psi_\ga(m)$ so that $c_\ga(s) = 2$. Hence
$\psi_\ga(m-1) = \psi_\gb(k-1)$. Also by the coding we have that
$c'_\ga(m)$ is the number of 0's in $c_\ga$ between $\psi_\ga(m-1)$ and
$\psi_\ga(m)$, which is the same as the number of 0's in $c_\gb$
between $\psi_\gb(k-1)$ and $\psi_\gb(k)$.

Finally, we can apply Lemma~\ref{tree} to get a strongly tree-like $\eta$ 
which satisfies $H$-uniformization.  \fin

\begin{lemma} 
\label{key}
Suppose that there is a ladder system $\zeta = \langle \zeta _\delta
\colon \delta \in \lim(\omega _1)\rangle $ such that $(\zeta,
2)$-uniformization holds, and suppose we are given a prime $p_\alpha$
for each $\alpha \in [\omega, \omega _1)$. Let $\{ x_\gn \colon \gn
\in \goo \}$ and $\{ y_\gn \colon \gn
\in \goo \}$ be sets of symbols.

 Then there are primes
$q_{\alpha ,n}$ for each $\alpha \in [\omega, \omega _1)$ and $n \in
\omega $ and a ladder system $\eta = \langle \eta _\alpha \colon
\alpha \in [\omega, \omega _1)\rangle$ such that given integers
$r_\alpha$ and $t_{\alpha ,n}$ for all $\alpha \in [\omega, \omega
_1)$ and $n \in \omega $, there is a function
$$\psi \colon \{x_\nu ,
y_\nu \colon \nu \in \omega _1\} \rightarrow {\open Z}$$
 such that for
all $\alpha \in [\omega, \omega _1)$ and all $n \in \omega $,
 $$\psi
(x_\alpha ) \equiv r_\alpha  \pmod{p_\alpha} {\rm \ \ and}$$
$$\psi (x_\alpha ) -
\psi (y_{\eta _\alpha (n)}) \equiv t_{\alpha ,n} \pmod{q_{\alpha
,n}}.$$
 Also, $\eta $ has the property that if $\eta _\alpha (m) = \eta
_\beta (n)$, then $m = n$, $p_\alpha = p_\beta $ and $\langle
q_{\alpha ,k}\colon k \leq n\rangle = \langle q_{\beta ,k}\colon k
\leq  n\rangle$. 
\end{lemma} 

\proof
Define the $q_{\alpha ,n}$ so that there is no repetition in the
sequence $\langle p_\alpha \rangle  \frown \langle q_{\alpha ,n}\colon
n \in \omega \rangle $ and such that if $p_\alpha  = p_\beta $, then
$q_{\alpha ,n} = q_{\beta ,n}$ for all $n$.  Without loss of
generality we can suppose that $r_\alpha  \in  \{0, \ldots, p_\alpha
-1\}$ and $t_{\alpha ,n} \in  \{0, \ldots, q_{\alpha ,n}-1\}$.  Fix a
bijection $\theta \colon {}^{<\omega }\omega  \rightarrow \omega$ such
that if $u\colon m \rightarrow  \omega $ and $v\colon m \rightarrow
\omega $ are such that $u(i) \leq v(i)$ for all $i < m$, then $\theta
(u) \leq  \theta (v)$. For each $\alpha $ and $n$, let $h_\alpha (n) =
\theta (\langle p_\alpha \rangle \frown \langle q_{\alpha ,j}\colon j \leq
n\rangle  )$
 Let $F$ be the function on $[\omega, \omega _1)$
such that $F(\alpha ) = p_\alpha $. Apply Lemma~\ref{prep} to this
situation to obtain the ladder system $\eta$ as in that lemma. Then
there is a uniformization $\langle f, f^*\rangle$ for the coloring given by
$c_\alpha (n) = \gth(\langle r_\alpha \rangle \frown \langle t_{\alpha ,j}\colon j
\leq  n\rangle $). 

	We can assume that $f^*(\ga)$ is minimal for $f$, i.e.,
$f^*(\ga)$ is the least $k$ so that $f(n) = c_\ga(n)$, for all $n
\geq k$. An immediate consequence of the minimality is that if
there exists $n
\geq f^*(\ga), f^*(\gb)$ with $\eta_\ga(n) = \eta_\gb(n)$ then
$f^*(\ga) = f^*(\gb)$. (The point is that $c_\ga(n) = c_\gb(n)$
implies that $c_\ga\rest n = c_\gb\rest n$.) 

	We now define $\psi $ in $\omega$ stages. At  stage $k$, we
will define $\psi(x_\alpha)$ for all
$\alpha $ such that $f^*(\alpha ) = k$ and we will define $\psi(y_\nu)$ 
for all $\nu $ of the form $\eta _\gamma(k)$ or of the
form $\eta _\alpha (n)$ where $f^*(\alpha ) = k$ and $n > k$.
First of all, for each $\nu $ of the form $\eta _\gamma (k)$ for some
$\gamma $, let $\psi (y_\nu )$ be arbitrary, if it has not already
been defined at a previous stage.  [Note that if $\nu $ is of this
form then $k$, but not $\alpha $, is uniquely determined by the
tree-like property of $\eta$.] For each $\alpha $ such that
$f^*(\alpha ) = k$, define $\psi (x_\alpha )$ to be the minimal
natural number such that 
$$ \psi (x_\alpha ) \equiv r_\alpha \pmod{p_\alpha}$$ 
$$ \psi (x_\alpha ) \equiv t_{\alpha ,j} + \psi
(y_{\eta _\alpha (j)}) \pmod{q_{\alpha ,j}}\hbox{ for }\ j
\leq  k.$$ 
This is possible by the Chinese Remainder Theorem. Now for
each $\nu $ such that $\nu  = \eta _\alpha (n)$ with $n > f^*(\alpha )
= k$, choose $\psi (y_\nu )$ minimal in $\omega $ such that $\psi
(x_\alpha ) - \psi (y_\nu ) \equiv  t_{\alpha ,n} \pmod{q_{\alpha
,n}}$; this is well-defined (independent of $\alpha $) by the
tree-like properties of $\eta $  and the primes, the uniformization,
and the minimal choices of $\psi (x_\alpha )$ and $\psi (y_\nu)$. 
Notice as well that by the minimality of $f^*$ and the remark above,
any $\nu$ is considered at at most one stage. To finish we let
$\psi(y_\nu)$ be arbitrary if $\nu$ is not of the form $\eta_\ga(n)$
for any $\ga$ or $n$. \fin

\begin{theorem} 
\label{main}
If there is a $W$-group A of cardinality $\aleph _1$ with $\Gamma (A)
= 1$, then there is a $W$-group $G$ of cardinality $\aleph _1$ which
is not strongly $\aleph _1$-free. 
\end{theorem} 

\proof
By Lemma~\ref{XII.3}, there is a ladder system $\zeta $ on $\lim(\omega
_1)$ which satisfies 2-uniformization. So we are in a position to
appeal to Lemma~\ref{key}. In fact by successive uses of this lemma,
we can define, by induction on $m \in  \omega$, sequences of primes
$\langle p^m_\alpha \colon \omega  \leq  \alpha < \omega _1\rangle $
and $\langle q^m_{\alpha ,n}\colon \omega  \leq  \alpha  <  \omega _1, n
\in \omega \rangle $, and ladder systems $\eta ^m = \langle \eta
^m_\alpha \colon \alpha  \in  [\omega, \omega _1)\rangle $ which for each
$m \in  \omega $ satisfy the properties given in Lemma~\ref{key} and
moreover are such that for all $m$, $\alpha $ and $n$, $p^{m+1}_{\eta
^m_\alpha (n)} = q^m_{\alpha ,n}$.  

Let $F$ be the free group on $\{x^m_\alpha \colon \omega  \leq  \alpha
< \omega _1, m \in  \omega \} \cup  \{z_{\alpha ,m,n}\colon \omega
\leq  \alpha <  \omega_1, m, n \in  \omega \}$ and let $K$ be
the subgroup of $F$ generated by $\{w_{\alpha, m,  n}\colon \omega  \leq
\alpha  < \omega _1, m, n \in  \omega \}$ where
$$
w_{\alpha ,m,n} = -q^m_{\alpha ,n}z_{\alpha ,m,n} + x^m_\alpha  - 
x^{m+1}_{\eta^{m} _\alpha (n)}.
$$
Let $G$ be $F/K$. In a harmless abuse of notation we shall identify
elements of $F$ with their images in $F/K = G$. To see that $G$ is not
strongly $\aleph _1$-free, consider the set $Y = \{x^m_\alpha \colon m <
\omega, \alpha  < \omega \} \subseteq  G$ and show by
induction on $\alpha < \go_1$ that if $H$ is an $\aleph_1$-pure subgroup of
$G$ containing $Y$, then $x^m_\alpha  \in  H$ for all $m \in  \omega
$. (The key point is that $x^m_\alpha $ will be divisible by infinitely
many primes modulo $H$ since $x^{m+1}_{\eta^m _\alpha (n)} \in  H$ by
induction.) 

To see that $G$ is a $W$-group, consider $f \in  \hom(K, {\open
Z})$. We want to define an extension of $f$ to $g \in  \hom(F,
{\open Z})$. The definition of $g$ will take place in $\omega$
stages. At the start of stage $k$,  for all $\alpha$ we have defined
$g(z_{\alpha ,m,n})$ for $m \leq  k-2$ and $g(x^m_\alpha )$ for $m
\leq  k - 1$, and  we have defined $r^k_\ga$ and committed  $g(x^k_\alpha )$
to be $r^k_\ga$ modulo $p^k_\alpha$  ($
= q^{k-1}_{\beta ,r}$ where $\eta ^{k-1}_\beta (r) =
\alpha$).

	Apply the uniformization property of Lemma~\ref{key} with
$r_\alpha = r^k_\ga$ and $t_{\alpha ,n} = f(w_{\alpha ,k,n})$.
We obtain a function $\psi_k \colon \{x^k_\nu $, $x^{k+1}_\nu \colon \nu
\in \omega _1\} \rightarrow {\open Z}$ such that $\psi_k (x^k_\alpha )
\equiv r^k_\ga \pmod{p^k_\alpha} $ and
$$\psi_k (x^k_\alpha ) -
\psi_k (x^{k+1}_{\eta^k _\alpha (n)}) \equiv f(w_{\alpha ,k,n}) \pmod{q^k_{\alpha ,n}}.$$
 Define $g(x^k_\alpha ) = \psi_k (x^k_\alpha )$ and
let $r^{k +1}_\ga$  be $\psi_k(x^{k+1}_\ga )$.
 Then by induction

\begin{eqnarray*}
g(x^{k-1}_\alpha) - g(x^k_{\eta^{k-1} _\alpha (n)}) & = &
 \psi_{k - 1}(x^{k - 1}_\ga) - \psi_k(x^k_{\eta_\ga^{k-1}(n)})\\
&\equiv &
 \psi_{k - 1}(x^{k - 1}_\ga) - r^k_{\eta^{k-1}_\ga(n)}
\equiv 
 \psi_{k - 1}(x^{k - 1}_\ga) - \psi_{k - 1}(x^k_{\eta_\ga^k(n)})\\
&\equiv &
  f(w_{\alpha ,k-1,n}) \pmod{q^{k-1}_{\alpha ,n}}.
\end{eqnarray*}

\noindent
 So define $g(z_{\alpha ,k-1,n})$ to be the
unique integer such that 
$$ g(x^{k-1}_\alpha ) - g(x^k_{\eta^{k-1} _\alpha (n)}) - f(w_{\alpha
,k-1,n}) = q^{k-1}_{\alpha ,n}g(z_{\alpha ,k-1,n}).  $$

\noindent
This completes the definition at stage $k$, and thus completes the proof.  
\fin

As mentioned before, Chase proved that CH implies that every
Whitehead group is strongly $\hao$-free. We can thus derive as a
consequence of Theorem~\ref{main} that CH implies that every Whitehead
 group $A$ of cardinality
$\hao$ satisfies $\gG(A) \neq 1$; this is a complicated way to prove
a fact already known, which is derived more easily using the weak
diamond principle (cf.\ \cite[XII.1.8]{EM}).

The following consequence of the theorem was also already known
 (see \cite[8.2, p.
74]{E}), but the proof here is more elegant, if less direct.

\begin{corollary}

There exists a Shelah group of cardinality $\ha_1$ which is not
strongly $\ha_1$-free.

\end{corollary}

\proof
Choose  sequences of primes
$\langle p^m_\alpha \colon \omega  \leq  \alpha < \omega _1\rangle $
and $\langle q^m_{\alpha ,n}\colon \omega  \leq  \alpha  <  \omega _1, n
\in \omega \rangle $, and ladder systems $\eta ^m = \langle \eta
^m_\alpha \colon \alpha  \in  [\omega, \omega _1)\rangle $
satisfying all the conditions in Theorem~\ref{main} except for the
uniformization properties. This can clearly be done in ZFC. Construct
$G$ as in Theorem~\ref{main}. Then, as before, $G$ is not strongly
$\ha_1$-free. We need to show that $G$ is a Shelah group. Note that
the property of {\it not} being a Shelah group of cardinality $\ha_1$
is absolute for extensions 
which preserve $\ha_1$.
 There is a generic extension of the universe
which satisfies MA + $\neg$CH. In this model, every ladder system
satisfies $\ha_0$-uniformization (cf.\ \cite[VI.4.6]{EM}), so our
ladder systems have the property given in Lemma~\ref{key}. Then the
proof of Theorem~\ref{main} applies to show that $G$ is a W-group. But in a
model of MA + $\neg$CH, every W-group is a Shelah group (cf.\
\cite[XII.3.20]{EM}). So $G$ was a Shelah group to begin with.
\footnote{referee: ``Remark that being `Shelah' is absolute''. {\em I have
not made any change here since there is already a remark.}}
  \fin

Combining Theorem~\ref{main} with results from \cite[Chapter
XII]{EM} we have a proof of Theorem~\ref{sumup} stated in the Introduction.

In a similar way one can also prove

\begin{theorem}
\label{sumup2}
The following are equivalent:

\noindent
(a) There is an $\aleph _1$-separable $\ha_1$-coseparable group A of cardinality $\aleph _1$ with 
$\Gamma (A) = 1$.

\noindent
(b) There is an strongly $\aleph _1$-free $\ha_1$-coseparable group A of cardinality $\aleph _1$ 
with $\Gamma (A) = 1$.

\noindent
(c) There is an $\ha_1$-coseparable group A of cardinality $\aleph _1$ with $\Gamma (A) = 1$.

\noindent
(d) There is an  $\ha_1$-coseparable group of cardinality 
$\aleph _1$ which is {\em not} strongly $\aleph _1$-free.

\noindent
(e) There is a ladder system on a stationary subset of $\ \lim (\omega _1)$ 
which satisfies $\omega $-uniformization. 

\end{theorem}





\section{The Second Question}

It is consistent that there are non-free Whitehead groups of
cardinality $\aleph _1$ but every $\aleph _1$-coseparable
group of cardinality  $\hao$ is free (see \cite[XII.3.18]{EM}). 
 Here we shall show that if {\em every}
strongly $\aleph _1$-free group of cardinality $\aleph _1$ is
Whitehead, then every ladder system on a stationary subset of $\lim
(\omega _1)$ has $\omega $-uniformization, and hence it follows that
there are non-free $\aleph _1$-coseparable groups of cardinality
$\aleph _1$.

\begin{proposition}

\label{A1} 

Assume that every strongly $\aleph _1$-free group of cardinality
$\aleph _1$ is Whitehead. Then for any ladder system $\eta  = \langle
\eta _\delta  \colon  \delta \in  S\rangle $ on a stationary subset
$S$ of $\lim
(\omega _1)$, and any $\omega $-coloring $c = \langle c_\delta  \colon 
\delta  \in  S\rangle $ of $\eta $, there is a pair $\langle g, g^*\rangle$ such
that $g^*\colon  S \rightarrow  \omega $ and $g\colon  \omega _1 \rightarrow
\omega $ such that for all $\delta  \in  S$ and all $n \in  \omega $,
if $n \geq  g^*(\delta )$, then $g(\eta _\delta (n)) > c_\delta
(n).$

\end{proposition} 

\proof Given what we are trying to prove, we can assume that each
$c_\delta $ is a strictly increasing function: $\omega  \rightarrow
\omega $. For each $\delta, n$ choose a prime $p_{\delta ,n} >
4c_\delta (n)$. Define $G$ to be the free group on $\{y_{\delta ,n}\colon 
\delta  \in  S$, $n \in  \omega \} \cup \{x_\nu \colon  \nu  \in  \omega
_1\}$ modulo the relations 
\begin{equation} p_{\delta ,n}y_{\delta
,n+1} = y_{\delta ,0} + x_{\eta _\delta (n)}.  
\label{G}
\end{equation}
It is routine to check that $G$ is strongly $\aleph
_1$-free. Let $H$ be the free group on $\{y'_{\delta ,n}\colon  \delta  \in
S$, $n \in  \omega \} \cup \{x'_\nu \colon  \nu  \in  \omega _1\} \cup
\{z\}$ modulo the relations 
\begin{equation} p_{\delta ,n}y'_{\delta
,n+1} = y'_{\delta ,0} + x'_{\eta _{\delta}(n)} + c_\delta (n)z.
\label{H} 
\end{equation}
Then there is a homomorphism $\pi $ of $H$ onto $G$ taking
$y'_{\delta,n}$ to $y_{\delta ,n}$ and $x'_\nu $ to $x_\nu $ and
which has kernel ${\open Z}z$.  By hypothesis, since $G$ is Whitehead,
there is a splitting $\varphi \colon  G \rightarrow  H$, i.e., such that $\pi
\circ  \varphi  = 1_G$. In particular, for all $\alpha  \in  \omega
_1$, there is $d(\alpha ) \in  {\open Z}$ such that $\varphi (x_\alpha
) - x'_\alpha  = d(\alpha )z.$ 

Define $g(\alpha ) = 2|d(\alpha )|$. Applying $\varphi $ to
equation~(\ref{G}) and subtracting (\ref{H}), we see that $p_{\delta
,n}$ divides 
\begin{equation} \varphi (y_{\delta ,0}) - y'_{\delta ,0}
+ \varphi (x_{\eta _\delta (n)}) - x'_{\eta _\delta (n)} - c_\delta
(n)z
\label{result} 
\end{equation}
 in $ {\open Z}z$.  Let $b$ be such that $bz = \varphi (y_{\delta ,0}) -
y'_{\delta ,0}$. Define $g^*(\delta )$ so that $c_\delta (g^*(\delta ))
> 2b$. 

	Assume that $n \geq g^*(\gd)$.  Then $p_{\delta ,n} >
 4c_\delta (n) > 8b.$ 
  Now consider two cases. The first
is that (\ref{result}) is zero, in which case $d(\eta _\delta (n))z =
\varphi (x_{\eta _\delta (n)}) - x'_{\eta _\delta (n)} = c_\delta (n)z
- bz$. Since $c_\delta (n) - b > c_\delta (n) - c_\delta (n)/2 =
c_\delta (n)/2$,  $c_\delta (n) < 2d(\eta _\delta (n))$, and thus
$c_\delta (n) < g(\eta_\delta(n)) $
 In the second case, (\ref{result})
equals $mz$ where $m$ is at least $p_{\delta ,n}$ in absolute
value, so $|d(\eta _\delta (n))| + |b - c_\delta (n)| \geq  p_{\delta
,n}$. But $|b - c_\delta (n)| \leq  |b| + |c_\delta (n)| < p_{\delta
,n}/2$, so $|d(\eta _\delta (n))| \geq p_{\delta ,n}/2 > c_\delta
(n)$. Hence $c_\delta (n) < g(\eta_\delta(n))$.  \fin 

\begin{corollary}
\label{corA1}

Assume that every strongly $\aleph _1$-free group of cardinality
$\aleph _1$ is Whitehead. Given a ladder system $\eta  = \langle
\eta _\delta  \colon  \delta \in  S\rangle $ on a stationary
 subset $S$ of $\lim
(\omega _1)$, there is a function
 $g\colon  \omega _1 \rightarrow
\omega $ such that for all $\delta  \in  S,$
 $g(\eta _\delta (n)) \geq n$ for all but finitely many $n \in \go$.
\end{corollary}

\proof
  Define an $\go$-coloring
 $c = \langle c_\gd\colon \gd \in S\rangle$ by $c_\gd(n) =
n - 1.$  There is a pair $\langle g, g^*\rangle$ as in
 Proposition~\ref{A1} with respect to
$c$. Clearly $g$ is the desired function.  \fin

\begin{lemma} 

\label{numbrthy} 

Given any positive integer $k$ and prime $p > 8k$, there are integers
 $a_0$ and $a_1$ and a function $F\colon 
{\open Z}/p{\open Z} \rightarrow  2$ such that for all $m \in  {\open
Z}$, if $|m| \leq  k$, then $F((m + a_\ell ) + p{\open Z}) = \ell $
for $\ell  = 0, 1$. 

\end{lemma} 

\proof Let $a_0 = 0$, $a_1 = 3k$. Then $\{m + a_0\colon  |m| \leq  k\} =
[-k$, $k]$ and $\{m + a_1\colon  |m| \leq  k\} = [2k$, 4k]. Since $p > 8k$,
$\{i + p{\open Z}\colon  -k \leq i \leq  k\}$ is disjoint from $\{j +
p{\open Z}\colon  2k \leq  j \leq  4k\}$, so we can define $F$ as desired.
\fin 

As mentioned in the Introduction, it was shown in \cite{Sh98} that
 if there is one
 strongly $\aleph
_1$-free group of cardinality $\aleph _1$ which is not free but
Whitehead, then there is some ladder system on a stationary subset of
$\omega _1$ which satisfies 2-uniformization.
Here we show: 

\begin{theorem} 

\label{main4A}

Assume that every strongly $\aleph _1$-free group of cardinality
$\aleph _1$ is Whitehead. Then every ladder system on a stationary
subset of $\lim (\omega _1)$ satisfies $2$-uniformization.  

\end{theorem} 

\proof Given a ladder system $\eta  = \langle \eta
_\delta  \colon  \delta \in  S\rangle $, let
$g$ be as in Corollary~\ref{corA1}.
By omitting a finite initial segment of each ladder,
 we can assume, without loss of
generality, that  $g(\eta
_\delta (n)) \geq n$ for all $n \in \go$.

For each $\alpha  \in  \omega _1$, choose a prime $p_\alpha > 8
g(\ga)$.
Also, for each $\alpha  \in  \aleph _1$, choose a function 
$$
F_\alpha \colon  {\open Z}/p_\alpha {\open Z} \rightarrow  2
$$
and integers $a^\alpha _0$, $a^\alpha _1$ such that for all
$m \in  {\open Z}$, if $|m| \leq  g(\alpha )$, then $F_\alpha (m +
a^\alpha _\ell ) = \ell $, for $\ell  = 0$, 1. (Here, and hereafter,
we  write $F_\alpha (k)$ instead of $F_\alpha (k + p_\alpha {\open
Z}).)$ This is possible by Lemma~\ref{numbrthy}. 

Now given a 2-coloring $c = \langle
c_\delta  \colon  \delta  \in S\rangle$ of $\eta$ define, as in
 Proposition~\ref{A1},  $G$ to be the free group on $\{y_{\delta ,n}\colon  \delta  
\in  S$, $n \in  \omega \} \cup  \{x_\nu \colon  \nu  \in  \omega _1\}$
 modulo the 
relations
\begin{equation}
p_{\eta _\delta (n)}y_{\delta ,n+1} = y_{\delta ,0} + x_{\eta _\delta (n)}.
\label{G2}
\end{equation}
and let $H$ be the free group on $\{y'_{\delta,n}\colon  \delta 
 \in  S$, $n \in  
\omega \} \cup  \{x'_\nu \colon  \nu  \in  \omega _1\} \cup \{z\}$ modulo the 
relations
\begin{equation}
p_{\eta _\delta (n)}y'_{\delta,n+1} = y'_{\delta,0} + 
x'_{\eta _\delta (n)} + a_{\delta ,n}z
\label{H2}
\end{equation}
where $a_{\delta ,n} = a^{\eta _\delta (n)}_{c_\delta (n)}$.
Let $\pi \colon H \longrightarrow G$ be the homomorphism taking
$y'_{\delta,n}$ to $y_{\delta ,n}$ and
$x'_\nu $ to $x_\nu $;
then there is a splitting $\varphi \colon  G \rightarrow
H$ of 
 $\pi $. We shall identify the elements of ${\open Z}z$
with integers; thus, for example, $\varphi (x_\alpha ) - x'_\alpha$
is an integer. 

Define the uniformizing function $f\colon  \aleph _1 \rightarrow
  2$ by
$$
f(\alpha ) = F_\alpha (\varphi (x_\alpha ) - x'_\alpha ). 
$$
We claim that $f(\eta _\delta (n)) = c_\delta (n)$ when $n \geq  
|\varphi (y_{\delta ,0}) - y'_{\delta,0}|$. As in Proposition~\ref{A1}, by 
applying $\varphi $ to (\ref{G2}) and subtracting (\ref{H2}), we get that  
$\varphi (x_{\eta _\delta (n)}) - x'_{\eta_\delta (n)}$ is congruent to 
$y'_{\delta,0} - 
\varphi (y_{\delta ,0}) + a_{\delta ,n}  \pmod{
p_{\eta _\delta }(n)}$. Hence
$$
f(\eta _\delta (n)) = F_{\eta _\delta (n)}(y'_{\delta,0} - 
\varphi (y_{\delta ,0}) + a_{\delta ,n})
$$
which equals $c_\delta (n)$ when $|y'_{\delta ,0} - \varphi
(y_{\delta ,0})| \leq  g(\eta _\delta (n))$ by choice of $F_{\eta
_\delta (n)}$. But in fact this is the case when $n \geq  |\varphi
(y_{\delta ,0}) - y'_{\delta,0}|$ because $g(\eta _\delta (n)) \geq n$.
\fin

\begin{lemma}
\label{disjtladders}

Given a stationary subset $S$ of $\lim(\go_1)$, for each $\ga \in
\go_1$ let $\gs(\ga)$ denote the least element of $S$ which is
greater than $\ga$. Then for each $\alpha  \in \omega _1$
there is a ladder $\zeta _\alpha $ on $\sigma
(\alpha )$ such that $\zeta _\alpha (0) > \alpha $ and such that for
all $\alpha  \neq  \beta $, $\rge(\zeta _\alpha ) \cap  \rge(\zeta
_\beta ) = \emptyset$.

\end{lemma}

\proof
For each $\grg \in S$, let $\grg^+$ denote the next largest element of
$S$. Then  $\gs(\ga) = \grg^+$ if and only if $\ga \in
 [\grg, \grg^+)$. It is clear that  $\grg^+$ 
contains the disjoint
union of $\go$  sets of order type $\go$, each of which is cofinal in
$\grg^+$:
\footnote{referee: `` $\amalg_{n , \go}$[in next line] should be
$\cup_{n < \go}$''. {\em I have left this alone.}}
$$\grg^+ \supseteq \amalg_{n < \go}W_n.$$

Let $\theta_\grg$ be a bijection of $[\grg, \grg^+)$ onto $\go$. Then
if $\ga \in [\grg, \grg^+)$,
let $\gz_\ga$ enumerate $W_{\theta_\grg(\ga)} \setminus (\ga + 1)$ in
increasing order.  \fin

The following result has been proved in \cite[1.\ 4, p.\ 262]{Sh98},
 but we 
give a self-contained proof here.

\begin{theorem} 

\label{ladders} 

Let $S$ be a stationary subset of $\lim (\omega _1)$. If every ladder
 system on
$S$ satisifes $2$-uniformization, then every ladder system on $S$
 satisfies 
$\omega $-uniformization. 

\end{theorem} 

\proof
\footnote{referee: `` A short intuitive remark will help the reader
understand the proof.''}
 Consider a ladder system $\eta  = \langle \eta _\delta \colon  \delta
\in  S\rangle $ and an $\omega $-coloring $c = \langle c_\delta \colon 
\delta  \in  S\rangle $.
We are going to define another ladder system $\eta' = \langle \eta
'_\delta \colon  \delta  \in  S\rangle $ and a 2-coloring $c'$.
Roughly, and slightly inaccurately, we get $\gh'$ from $\gh$ by adding a
segment of length $c_\gd(n)$ at each $\gh_\gd(n)$ and then we color the
new segment by a binary code for $c_\gd(n)$.  

Let the ladders $\zeta _\alpha $ be
as in Lemma~\ref{disjtladders}.
Let $\eta '_\delta $ enumerate the $\go$-sequence
$$ \cup _{n \in  \omega }\{\zeta _{\eta_\delta (n)}(k)\colon  k \leq
c_\delta (n)\}. $$  
Define 
$c'_\delta (k) = 0$ if $\gh_\gd'(k) =
 \zeta _{\eta _\delta (n)}(c_\gd(n))$ 
for some $n$, 
and $c'_\delta (k) = 1$
otherwise.

By hypothesis, there is a uniformization $\langle f, f^*\rangle$ of the coloring
$c'$ of $\eta'$. Define $g\colon  \omega _1 \rightarrow  \omega _1$ as
follows: $g(\alpha )$ equals the number of 1's before the first 0 in
$f\rest \rge(\zeta _\alpha )$.  Define $g^*\colon  S \rightarrow  \omega $
by: $g^*(\delta ) = m$ if $m$ is minimal such that for every $n <
f^*(\delta )$, there exists $k < m$ with $\eta '_\delta (n) \in  \zeta
_{\eta _\delta (k)}$.  

We claim that $\langle g, g^*\rangle$ uniformizes the coloring $c$ of $\eta $.
Suppose $m \geq g^*(\gd)$. Let $\langle n_i \colon i \leq
c_\gd(m)\rangle$ enumerate in increasing order the set
$$\set{j}{\eta'_\gd(j) \in
\rge(\zeta_{\eta_\gd(m)}{\rest}(c_\gd(m)+1))}.$$  Then $c'_\delta (n_i)
= f(\eta'_\delta (n_i))$ for $i \leq c_\gd(m)$. So there are exactly
$c_\gd(m)$ 1's before the first 0 in $f\rest
\rge(\zeta_{\eta_\gd(m)})$. \fin



We can now give the proof of Theorem~\ref{sumup3} stated in the Introduction:
(1) implies (2) is trivial; (2) implies (3) is Theorem~\ref{main4A};
(3) implies (4) is Theorem~\ref{ladders}; and (4) implies (5) is
a consequence of \cite[XII.3.1]{EM}.  \fin

 The third author claims to have a proof of 
(4) implies (1) and hence an
affirmative answer to the second question (in the Introduction); but
he has not yet been able to 
convince the first two authors.




\section{Uniformization on a cub}

The theorems of this section have no direct application to Whitehead
groups, but they complete a circle of results regarding uniformizations.

\begin{theorem} 
\label{cub}

Suppose that $S$ is a stationary subset of $\lim (\omega _1)$ which
has the property that for every ladder system
 $\eta  = \langle \eta _\delta \colon  \delta
\in  S\rangle $ on  $S$  and  every $\omega
$-coloring $c = \langle c_\delta \colon  \delta  \in S\rangle $, there is a
pair $(f, f^*)$  and a cub $C$ on $\omega _1$ such that for every
$\delta  \in  S \cap  C$, $f(\eta _\delta (n)) = c_\delta (n)$ for all
$n \geq  f^*(\delta )$. Then every ladder system on $S$  satisfies
$\omega $-uniformization.

\end{theorem} 

\proof Let $\eta $ be as given and let $c$ be any $\omega $-coloring
of $\eta $. Let $C$ and  $(f, f^*)$ be as in the statement of the
theorem. For each $\ga \in C$, let $\gth_\ga$ be a bijection from $\go$
onto $\ga$.\footnote{Alan: note addition(mine); somehow this sentence
got lost in later versions.}

Let $S_1 = C^* \cap S$, where $C^*$ is the
set of limit points of $C$. For each $\delta  \in S_1$, let $\eta
^1_\delta $ enumerate in increasing order the set $\cup _{n \in
\omega }Z_n$, where $Z_n$ is defined as follows.  Let $\gamma _n =
\min (C \setminus  (\eta _\delta (n) + 1))$, i.e., $\gamma _n$ is
 the
least element of $C$ which is greater than $ \eta_\delta (n)$;
 then
$\eta _\delta (n) = \theta _{\gamma _n}(k_n)$ for some unique $k_n \in
\omega $. Define
$$
Z_n = \{\sigma \colon   \sigma  = \theta _{\gamma _n}(j)
\hbox{ for some }\ j \leq  
k_n\hbox{ and }\ \min (C \setminus  (\sigma  + 1)) = \gamma _n\}
$$
Note that $\eta _\delta (n) \in  Z_n$, so the range of $\eta ^1_\delta
$ includes the range of  $\eta _\delta $. We are going to define a
coloring $c^1 = \langle c^1_\delta \colon  \delta  \in  S_1\rangle $. It
will be convenient to regard $c^1_\delta $ as a function whose domain
is $\rge(\eta ^1_\delta )$ rather than $\omega $; that is, if $\sigma
= \eta ^1_\delta (k)$, we shall write $c^1_\delta (\sigma )$ instead
of $c^1_\delta (k)$. For all $\delta  \in S_1$ and $n \in  \omega $,
if $\eta ^1_\delta (n) = \eta _\delta (m)$,
then $c^1_\delta
(\eta ^1_\delta (n))$ is defined to be
 $$\langle \langle \eta _\delta
(j), r_{\gd,j} \rangle
 \colon  j
\leq  m\rangle$$
where $r_{\gd, j}$ is the size of the intersection of the open interval
$(\eta_\gd(j), \min(C \setminus (\eta_\gd(j)+1)))$ with
$\rge(\eta_\gd)$;
$c^1_\delta (\eta ^1_\delta (n))$    
 can be regarded as
an element of $\omega  \setminus  \{0\}$ by a coding argument.
Otherwise $c^1_\delta (\eta ^1_\delta (n))$ is defined to be 0.

By hypothesis
\footnote{referee: `` Page 15, lines 11b -1b: It seems that you assume
some maximality property for $f_1$, something like: `if for any $\gd \in
D_1$ and $\grg \in \rge(\gh^1_\gd)$ $c^1_\gd(\grg)$ does not depend on
$\gd$, then $f_1(\grg)$ is this constant value.'
(Without this assumption $f_1$ may get `wild' values on $\grg$'s which
appear in $\gh^1_\gd$ with index less than $f^*_1(\gd)$).''
{\em Alan: I don't see the point since we use only $f_1\rest A \cup B$ .}}
 there is a pair $(f_1, f^*_1)$ and a cub $C_1$ such
that for $\delta  \in  S_1 \cap  C_1$ and $n \geq  f^*_1(\delta )$,
$f_1(\eta ^1_\delta (n)) = c^1_\delta (\eta ^1_\delta (n))$. Without
loss of generality, we can assume that $C_1 \subseteq  C^*$.  

Define $D_1 = C_1 \cap  S$, $D_2 = (C \cap  S) \setminus  C_1$, $D_3 =
S \setminus  (D_1 \cup  D_2)$. We are going to define the desired
uniformization, $(f_0$, $f^*_0)$, of $c$ by defining $f^*_0 = g^*_1
\cup  g^*_2 \cup  g^*_3$, where $g^*_i\colon  D_i \rightarrow  \omega $. 

Define $g^*_1(\delta )$ to be the maximum of $f^*(\delta )$ and the
least $m$ such that $\eta _\delta (m) \geq  \eta ^1_\delta(k)$, where
$k \geq f^*_1(\delta )$ and there is $\ga \in C$ so that
$\eta_\gd^1(k-1) \leq\ga < \eta_\gd^1(k)$.
 Thus if $n \geq  g^*_1(\delta )$, $f(\eta _\delta (n)) =
c_\delta (n)$ and $f_1(\eta _\delta (n)) = c^1_\delta (\eta _\delta
(n))$. Let 
$$ A = \{\eta _\delta (n)\colon   \delta  \in  D_1,  n
\geq  g^*_1(\delta )\}.  $$ 
Let $\alpha  < \beta $ be two successive
members of $C$ (so, in particular, $(\alpha, \beta ] \cap  C_1 =
\emptyset $, since  $C_1 \subseteq  C^*)$. Notice that,
\footnote{referee: `` Page 15 line 1b - Page 16 line 1:
Explain the statement `Notice that...' ''; {\em Alan:
note changes in this sentence}}
 by the last clause in the definition of
$g^*_1$, if for some $\gd$ and $n$, $\eta_\gd(n) \in A \cap (\ga, \gb]$, then
$c^1_\gd(\grg) = f_1(\grg)$ for all $\grg \in (\ga,
\gb]\cap \rge(\eta^1_\gd)$.
We claim that there exists $\delta  \in  S_1$ such that $\delta  >
\beta $ and $A \cap  (\alpha , \beta )$ is contained in $\rge(\eta
_\delta )$; this implies that  $A \cap  (\alpha, \beta )$ is finite.
It suffices to show that for any $\delta _1$, $\delta _2$ in $D_1
\setminus  (\beta  + 1)$, if $\rge(\eta _{\delta _\ell }) \cap   A
\cap  (\alpha $, $\beta ) \neq \emptyset $ for $\ell = 1, 2,$ then
$\rge(\eta _{\delta _1})\cap \beta  =\rge(\eta _{\delta _2})\cap
\beta$
 Now, for each $\ell $, $n$ such
that $\eta _{\delta _\ell }(n) \in  (\alpha, \beta )$, note that $\min
(C \setminus  (\eta _{\delta _\ell }(n) + 1)) = \beta $, so there is a
$k_{\ell ,n}$ such that $\theta _\beta (k_{\ell ,n}) = \eta _{\delta
_\ell }(n)$ and a set 
$$Z^\ell _n = \{\sigma \colon   \sigma  = \theta _\beta (j) \hbox{ for
some }\  j \leq  k_{\ell ,n}
 \hbox{ and }\ \min (C \setminus  
(\sigma  + 1)) = \beta \} \subseteq  \rge(\eta ^1_{\delta _\ell }).$$
The sets $Z^\ell_n$ are linearly ordered by 
inclusion, so for each $\ell$, there is a largest one, which we shall
denote $Z^\ell $. Without loss of generality, $Z^1 \subseteq  Z^2$. 

	Also by the choice of $g_1$, we know for each $\ell$ that
$f_1(\gs) = c_{\gd_\ell}^1(\gs)$ for $\gs \in Z^\ell$. So for any $\gs \in
\rge(\eta_{\gd_1}) \cap (\ga, \gb)$, $0 \neq c_{\gd_1}^1(\gs) =
f_1(\gs) = c_{\gd_2}^1(\gs)$. Hence $\rge(\eta_{\gd_1}) \cap \gb\se
\eta_{\gd_1} \cap \gb$. Finally if we choose $m$ maximal so that
$\eta_{\gd_1}(m) \in (\ga, \gb)$ then $|(\eta_{\gd_1}(m), \min(C
\setminus (\eta_{\gd_1}(m)+1)) \cap \rge(\eta_{\gd_1})|= 0$.
Then $c_{\gd_1}^1(\gh_{\gd_1}(m)) =
f_1(\gh_{\gd_1}(m)) = c_{\gd_2}^1(\gh_{\gd_1}(m))$, so by definition of
 $c^1$, $\gh_{\gd_1}(j) = \gh_{\gd_2}(j)$ for all $j \leq m$ and
$\gh_{\gd_2}(m)$ is the largest element of $\rge(\gh_{\gd_2}) \cap
 \gb.$
 So we are done.

Define $h^*_2\colon  D_2 \rightarrow  \omega $ such that if
 $\alpha  < \beta $ are 
successive members of $C_1 \cup  \{0\}$, and if we define
$$
B_{\alpha ,\beta } = \{\eta _\delta (n)\colon  \delta  \in  D_2 \cap  
(\alpha,  \beta ),  h^*_2(\delta ) \leq  n < \omega \},
$$
then $B_{\alpha,\beta } \subseteq  (\alpha,\beta )$ and for any 
$\delta  \in  D_3$, $\sup (B_{\alpha ,\beta } \cap  \delta ) < \delta $. This
is not hard to do. Now 
define $g^*_2(\delta ) = \max \{h^*_2(\delta )$, $f^*_1(\delta )\}$ for all 
$\delta  \in  D_2$. Let\footnote{referee: ``page 16 line 16b: $f^*_4$
should be $f^*_0$.'' {\em Alan: I don't know what he means; my Xerox of
the printed copy he has does not have an $f^*_4$ that I can see.}}
$$
B = \{\eta _\delta (n)\colon  \delta  \in  D_2, g^*_2(\delta ) \leq  n < 
\omega \}.
$$
Thus for any $\delta  \in  D_3$, $B \cap  \delta $ and $A \cap  \delta $ are 
bounded in $\delta $ (the latter because there are successive elements $\alpha 
< \beta $ in $C$ such that $\alpha  < \delta  < \beta $ --- since $\delta  
\notin  C)$. Define $g^*_3\colon  D_3 \rightarrow  \omega $ such that for all 
$\delta  \in  D_3$, $\{\eta _\delta (n)\colon  g^*_3(\delta )
 \leq  n < \omega  \} 
\cap  (A \cup  B) = \emptyset$. 

Then let $f^*_0 = g^*_1 \cup  g^*_2 \cup  g^*_3$. We can then let $f_0\rest A 
\cup  B = f_1\rest A \cup  B$ and easily define $f_0$ on $\omega _1 \setminus  
(A \cup  B)$ to take care of those $\delta $ in $D_3$. \fin

We shall abbreviate the property given in the hypothesis of Theorem~\ref{cub}
by saying ``every ladder system on $S$ satisfies
 $\go$-uniformization on a cub''.
Combining the results of this section with those of the previous
section we have the following.

\begin{theorem}
\label{sumup4}

Let $S$ be a stationary subset of $\lim(\go_1)$.
Consider the following hypotheses.

(1) Every strongly $\ha_1$-free group $A$ of cardinality $\ha_1$
with $\gG(A) \subseteq S$ is
 $\ha_1$-coseparable.

(2)  Every strongly $\ha_1$-free group of cardinality $\ha_1$
 with $\gG(A) \subseteq S$ is
Whitehead.

(3) Every ladder system on $S$ satisfies
$2$-uniformization.

(4) Every ladder system on $S$ satisfies
$\go$-uniformization.

(5) Every ladder system on $S$ satisfies
$2$-uniformization on a cub.

(6) Every ladder system on $S$ satisfies
$\go$-uniformization on a cub.

\noindent
Then (1) $\Rightarrow$ (2) $\Rightarrow$ (3) $\Leftrightarrow$
 (4) $\Leftrightarrow$
(5) $\Leftrightarrow$ (6).

\end{theorem}

\proof
(1) implies (2),  (4) implies (3), (6) implies (5)
and (4) implies (6) are trivial. 
 Inspection of the
proof of Theorem~\ref{main4A}
 shows that it ``localizes'' to $S$, so
 (2) implies (3).  The implication from (3) to (4) is
 Theorem~\ref{ladders}
 the proof of (5) implies (6) is
exactly the same. That (6) implies (4) is Theorem~\ref{cub}. \fin



\begin{quote}
\it in the file received from Eklof, this section 4 was commented out.
---Martin
\end{quote}

\section{Topological Considerations}

	Uniformization results have been associated with the
construction of interesting normal
spaces. From the existence of a ladder system with
2-uniformization it is easy to construct a
normal space which is not metrizable. In fact
this is how the consistency with GCH of the failure of
 the normal Moore space conjecture was established \cite{ShW1}.
(See~\cite{T} for more information about the normal Moore space
conjecture.) A key difference between the Whitehead problem and the
construction of normal spaces from ladder systems is that in the
topological case the proof of the normality of the space does not
require the full power of $2$-uniformization, but only requires
uniformization of monochrome colourings. However, by considering a
large collection of spaces built from ladder systems, we can get
topological equivalents to uniformization principles. We would like to
thank Frank Tall for looking at this section, saving
us from an elementary error in topology, and providing information
about the normal Moore space conjecture.

        Recall that if $\ga = \gd + n$ where $\gd \in \lim(\go_1)$,
then a ladder on $\ga$ is defined to be a ladder on $\gd$.  Suppose
that $E \se [\go, \go_1)$ and $\eta$ is a ladder system on $E$. Then
we define a topological space $X(\eta)$ on $\go_1$ by defining by
induction on $\ga < \go_1$ a neighborhood base of $\ga$.  Let $\ga$ be
isolated if $\ga \notin E$. If $\ga \in E$, then a neighbourhood base
of $\ga$ is formed by the sets $\{\ga\} \cup \bigcup_{n\leq m} u_m$
where $u_m$ is a neighbourhood of $\eta_\ga(m)$ and $n < \go$.

 Suppose that $S \se \lim(\go_1)$. Let
$K_0(S)$ be the set of topological spaces of the form $X(\eta)$ where $\eta$
is a ladder system on some $E \se \set{\gd +n}{\gd \in S \mbox{ and }
n \in \go}$ which satisfies the additional hypothesis that if $\gd
+n, \gd +m \in E$ and $m \neq n$ then $\rge(\eta_{\gd+n}) \cap
\rge(\eta_{\gd+m}) = \emptyset$. Let $K_1(S)$ be the subset of
$K_0(S)$ consisting of all $X(\eta)$ such that
if $\eta = \{\eta_\ga \colon \ga \in E\}$, then
 for all $\ga \in E$,
 the range of $\eta_\ga$ consists of isolated points (i.e., elements of
$\go_1 \setminus E$).

	These classes of spaces can be used to give equivalents to
uniformization principles.

\begin{theorem}
Let $S\se \lim(\go_1)$.  The following are equivalent.

\noindent
(a) every ladder system on $S$ satisfies $2$-uniformization;

\noindent
(b) every member of $K_0(S)$ is normal;

\noindent
(c) every member of $K_1(S)$ is normal;
    
\noindent
(d) every ladder system on $S$ satisfies $\ha_0$-uniformization.

\end{theorem}

	The equivalence of (a) and (d) has already been
established. The rest of the section is devoted to
proving the non-trivial implications.

        From now on we will assume that every ladder system on a set
$E \se [\go, \go_1)$ is such that if  $\gd
+n, \gd +m \in E$ and $m \neq n$ then $\rge(\eta_{\gd+n}) \cap
\rge(\eta_{\gd+m}) = \emptyset$. With this assumption, there is a
simple connection between uniformization on subsets of $\lim(\go_1)$
and subsets of $[\go, \go_1)$. 

\begin{proposition}
    Suppose $S\se \lim(\go_1)$ and every ladder system on $S$
satisfies $2$-uniformization ($\ha_0$-uniformization).  If $E
\se \set{\gd+n}{\gd \in S \mbox{ and }n \in \go}$, then every ladder
system on $E$ satisfies $2$-uniformization ($\ha_0$-uniformization).
\end{proposition}

\proof Given $\set{\eta_\ga}{\ga \in E}$, for each $\gd \in S$
choose $\eta_\gd^*$ so that for all $n$, if $\gd +n \in E$ then the
range of  $\eta_{\gd+n}$ is contained, except for a finite set, in the
range of $\eta_\gd^*$. Let $\eta^* = \set{\eta_\gd^*}{\gd \in S}$.
Given a coloring $c =  \set{c_\ga}{\ga \in E}$ it is easy to produce a
colouring $c^*$ of $\eta^*$ such that any function which uniformizes $c^*$
also uniformizes $c$. \fin

        If $S$ is a stationary subset of $\lim(\go_1)$ and $\eta$ is a
ladder system on $S$ such that the ladders consist of successor
ordinals, then the space $X(\eta)$ is not metrizable.

The connection with the normal Moore space problem came from the
following easy fact.

\begin{theorem}
    Suppose $E\se [\go, \go_1)$  and $\eta$ is a ladder system on $E$
which satisfies $2$-uniformization where for all $\ga
\in E$ the range of $\eta_\ga$ consists of
isolated points. Then the space $X(\eta)$ is normal.
\end{theorem}

\proof
 Suppose $A_0$ and $A_1$ are disjoint closed sets. Choose
a coloring so that $c_\ga$ is constantly $0$ if $\ga \in A_0$ and
$c_\ga$ is constantly $1$ if $\ga \in A_1$. Suppose that $f$
uniformizes the coloring. Then we can let $U_0 = A_0 \cup
\set{\gb}{f(\gb) = 0 \mbox{ and }\gb \notin (E \cup A_1)}$ and $U_1 = A_1 \cup
\set{\gb}{f(\gb) = 1 \mbox{ and }\gb \notin (E \cup
A_0)}$.
 \fin

	Unlike the case of abelian groups, where the group
constructed from the ladder system is a Whitehead group if and only
if the ladder system has $2$-uniformization, we cannot deduce the
converse here because in the topological case we only need to deal with
monochromatic colorings.

\begin{theorem}
    Suppose $S \se \lim(\go_1)$. If every element of $K_1(S)$ is
normal then every ladder system on $S$ satisfies $2$-uniformization.
\end{theorem}

\proof Suppose we are given $\eta = \set{\eta_\gd}{\gd \in S}$ a
ladder system on $S$ and $c = \set{c_\gd}{\gd\in S}$ a coloring of
$\eta$. Let $\set{\gz(\ga)}{\ga < \go_1}$ enumerate the ordinals
equivalent to $2 \pmod{3}$ in increasing order.
 For $\gd \in S$ and $i
= 0, 1$ if there exists infinitely many $n$ so that $c_\gd(n) = i$,
let $\eta_{\gd+i}^*$ enumerate $\set{\gz(\eta_\gd(n))}{c_\gd(n) = i}$
in increasing order. Otherwise $\eta_{\gd+i}^*$ is undefined. Let
$\eta^* = \set{\eta_{\gd+i}^*}{\eta_{\gd+i}^* \mbox{ is defined}}$. It
is easy to see that $X(\eta^*) \in K_1(S)$ and so by hypothesis is
normal. For $i = 0, 1$, let $A_i = \set{\gd+i}{\gd \in S \mbox{ and }
\eta_{\gd+i}^* \mbox{ is defined}}$. Let $U_i$ be as guaranteed by
normality and choose $f$ so that $f(\ga) = i$ if $\gz(\ga) \in U_i$.
It is easy to check that $f$ uniformizes $c$. \fin

The previous two results show that (a) is equivalent to (c). It remains
to prove that (a) implies (b).

\begin{theorem}
    Suppose $S \se [\go,\go_1)$ and every ladder system on $S$
satisfies $2$-uniformization. Then every element of $K_0(S)$ is
normal. 
\end{theorem}

\proof
 Fix $\eta = \set{\eta_\gd}{\gd \in S}$. It suffices to show that
if
$C, D$ are disjoint closed sets, then there exists
$C'$ so that: (0) $C \se C'$; (1)  $C'
\cap D = \emptyset$; (2) $C'$ is closed; and (3) for all $\ga
\in C$ there is $k$ so that $\eta_\ga(r) \in C'$ for all $r \geq
k$. Before proving that $C'$ exists let us see why the claim
suffices. 

         Given disjoint closed sets $A_0$ and $A_1$,
let $ A_{n0} = A_n$ for $n \in \{0, 1\}$. Considering each $n\in \{0, 1\}$
alternately, we can inductively define $A_{nm}$, such that if $n$ is the
number considered at stage $m$, then $A_{nm+1}$ is to $A_{nm},
A_{(1-n)m}$ as $C'$ is to $C, D$. We also let $A_{(1-n)m+1} =
A_{(1-n)m}$.
Then let $A^n = \bigcup_{m < \go} A_{nm}$. To finish the proof we
must show that $A^n$ is open. We do this by induction on $\ga \in
A^n$. Suppose $\ga \in A^n$ and choose a stage $m$ where $n$ is
considered and $\ga \in A_{nm}$. If $\ga$ is isolated then we are
done; otherwise $\eta_\ga$ is defined. So for some $k$ and all $r \geq
k$, $\eta_\ga(r) \in A_{nm+1} \se A^n$. By induction, $A^n$ contains
an open neighborhood $u_r$ of each $\eta_\ga(r)$. Hence $A^n$
contains $\{\ga\}\cup \bigcup_{k \geq r} u_r$, which is an open
neighborhood of $\ga$. 

        It remains to show that $C'$ exists. For $\ga \in S$,
define $c_\ga$ to be constantly 0 if $\ga \in C$ and let $c_\ga$ be
constantly 1 if $\ga \notin C$.  Choose $f$ which uniformizes the
coloring. Let $C' = C \cup\set{\ga}{f(\ga) = 0 \mbox{ and } \ga
\notin  D}$.  Requirements (0) and (1) follow from
the definition. For clause (2) we must show that the complement of
$C'$ is open. By induction we show that if $\gb \notin C'$ then
the complement of $C'$ contains an open neighborhood of $\gb$. If
$\gb$ is isolated, there is nothing to prove. Otherwise $\eta_\gb$
exists and $\gb \notin C$. Since $f$ uniformizes the coloring there
is $n_0$ so that for all $m \geq n_0$, $f(\eta_\gb(m)) =1$.
Furthermore since $C$ is closed there is $n_1$ so that for all $m
\geq n_1$, $\eta_\gb(m) \notin C$. So if we let $n = \max\{n_0,
n_1\}$, for all $m \geq n$, $\eta_\gb(m) \notin C'$.  By the
induction hypothesis there is an open neighborhood of each
$\eta_\gb(m)$ contained in the complement of $C'$. So the
complement of $C'$ is open. The verification of (3) is similar to
the verification of (2) except we use that $D$ is
closed as well as that $f$ uniformizes the coloring. \fin

\begin{quote}
\it In the file received from Eklof, there was an ``end document''
here. I guess this means the rest of this {\bf section} is not inteded
to be printed.  I print it anyway. ---Martin
\end{quote}

proof of fact that $X(\eta)$ is not metrizable:
Assume the space is metrizable with metric $d$. For
each $\gd$ in $S$, let $O_\gd$ be the open set 
consisting of the ladder and notice that if $\tau \in S$ then $\tau
\notin O_\gd$. Choose $\gre_\gd$ so that $B_{\gre_\gd} \se O_\gd$. For
all $\gd \in S$ choose $n$ so that $d(\gd, \eta_\gd(n)) < \gre_\gd/2$.
Define $f(\gd) = \eta_\gd(n)$. By the pressing down lemma there is
$\tau, \gd \in S$ so that $f(\tau) = f(\gd)$. We can assume that
$\gre_\tau \leq \gre_\gd$. So $d(\gd, \tau) \leq d(\gd, f(\gd)) +
d(f(\gd), \tau) < \gre_\gd/2 + \gre_\tau/2
\leq \gre_\gd$. Hence $\tau \in B_{\gre_\gd} \se O_\gd$, a
contradiction.

\end{document}